\documentclass[12pt]{article}
\usepackage{latexsym}
\title{On the growth of restricted integer partition functions}
\author{E. Rodney Canfield\\
University of Georgia\\
Athens, GA 30602-7404\and
Herbert S. Wilf\\
University of Pennsylvania\\
Philadelphia, PA 19104-6395}
\newtheorem{theorem}{Theorem}
\newtheorem{lemma}{Lemma}
\newcommand{\eqdef}{\, =\kern -12.7pt\raise 6pt\hbox{{\tiny\textrm{def}}}\,\,}
\begin{document}
\maketitle

\begin{abstract}
We study the rate of growth of $p(n,S,M)$, the number of partitions of $n$ whose parts all belong to $S$ and whose multiplicities all belong to $M$, where $S$ (resp. $M$) are given infinite sets of positive (resp. nonnegative) integers. We show that if $M$ is all nonnegative integers then $p(n,S,M)$ cannot be of only polynomial growth, and that no sharper statement can be made. We ask: if $p(n,S,M)>0$ for all large enough $n$, can $p(n,S,M)$ be of polynomial growth in $n$?
\end{abstract}

\medskip

\noindent\textbf{Subject classification number}: 05A17 Primary

\smallskip


\medskip

\section{The question}
Let $S$ be a set of positive integers, and let $p_S(n)$ denote the number of partitions of the
 integer $n$ all of whose parts lie in $S$. For various sets $S$, the asymptotic growth  rate
 of $p_S(n)$ is known, and the known rates lie in the range of polynomial growth to
 superpolynomial-but-subexponential rates.

For example, if $S$ consists of all positive integers then the celebrated theorem of Hardy,
 Ramanujan and Rademacher \cite{hr, rad} has given the complete asymptotic expansion, of which
 the first term is
\begin{equation}
\label{eq:hrr}
p_S(n)\sim \frac{1}{4n\sqrt{3}}\exp{\left(\pi \sqrt{\frac{2n}{3}}\right)}.
\end{equation}

 As an example of a sparse set of parts, take $S=\{1,2,2^2,2^3,\dots\}$, the case of
 binary partitions. Then de Bruijn \cite{deb} found several terms of the asymptotics of the
  logarithm, which begins as
\begin{equation}
\label{eq:deb}
\log{p_S(2n)}=\frac{1}{2\log{2}}\left(\log{\frac{n}{\log{n}}}\right)^2(1+o(1)).
\end{equation}

For a final example, suppose the set $S$ of allowable parts is finite. Then we have,  say,
$S=\{a_1<\dots <a_k\}$, and we are dealing with ``the money changing problem,'' a.k.a.
``the problem of Frobenius.'' A result of Schur \cite{sch} holds that in this case $p_S(n)$
is of polynomial growth.
\begin{theorem} [Schur]\label{th:schur} If  $S=\{a_1<\dots <a_k\}$, and gcd$(S)=1$, then
\begin{equation}
\label{eq:mcp}
p_S(n)\sim \frac{n^{k-1}}{(k-1)!a_1a_2\dots a_k},
\end{equation}
and in particular, $p_S(n)>0$ for all large enough $n$.
\end{theorem}

We show here that if the set of allowable parts is infinite, no
matter how sparse, then the partition function
 $p_S(n)$ must grow faster than every polynomial. We show also that this result is best
  possible in the sense that if $\epsilon(n)$ is any unbounded function of $n$ then
   there exists an infinite set $S$ of allowable parts such that
   $p_S(n)=O(n^{\epsilon(n)})$.

   We discuss also the situation in which we have an arbitrary set
   of allowable parts and an arbitrary set of allowable
   multiplicities.
\section{Preliminaries}
\begin{lemma}\label{lem:cop} Let $S=a_1<a_2<a_3<\dots$ be a set of positive integers such that
gcd$(S)=1$. Then $S$ contains a \textit{finite}
 coprime subset.
\end{lemma}
\noindent Proof. Let $g_n=$gcd$(a_1,\dots,a_n)$. Then
  $a_1\ge g_1\ge g_2\ge \dots$, so $\exists i_0$ such that
 $\forall i>i_0$: $g_i=1$. Indeed, if not then $\exists i_0$ such that
 $\forall i>i_0$: $g_i=g>1$. But then we would have
  gcd$(S)=g>1$, a contradiction. $\Box$

\begin{lemma}
\label{lem:ps}
The following two properties of a set $S$ of positive integers are equivalent:
\begin{enumerate}
\item \label{it1} for all sufficiently large integers $n$ we have $p_S(n)>0$
\item gcd$(S)=1$.
\end{enumerate}
\end{lemma}
\noindent Proof. If gcd$(S)=1$ then by Lemma \ref{lem:cop} $S$
contains a \textit{finite} coprime subset $\overline{S}$. By Schur's
theorem, $p_{\overline{S}}(n)>0$ for all large enough $n$, hence so
is $p_S(n)$, and conclusion \ref{it1} holds. On the other hand, if
gcd$(S)>1$ then
 conclusion \ref{it1} is obviously false. $\Box$

We remark that Lemma \ref{lem:ps}, whose proof we have given in order to keep this paper self-contained,  is a special case of a much more general result of Bateman and Erd\H os \cite{ber}, who found the conditions on $S$ under which, for a fixed $k\ge 0$, almost all values of the $k$th differences of $\{p_S(n)\}_{n=0}^{\infty}$ are strictly positive.

Next we will need a lemma that allows us to estimate the growth of $p_S(n)$ for
arbitrary sets $S$ of parts. We will in fact prove a more general result, in which
not only the set $S$ of allowable parts can be arbitrarily prescribed, but so can
the set $M$ of allowable multiplicities of those parts.

Hence, let $S$ be a set of positive integers and let $M$ be a set of nonnegative
integers such that $0\in M$. Let $M(x),S(x)$ denote the respective counting functions
 of $M,S$. That
is $M(x)=|\{\mu\in M:\mu\le x\}|,$ and likewise for $S(x)$. Finally we denote by
 $p(n;S,M)$ the number of partitions of $n$ whose parts all belong to $S$ and the
 multiplicities of whose parts all belong to $M$.

\begin{lemma}\label{lem:ineq} For the general partition function $p(n;S,M)$ we have
\begin{equation}
\label{eq:upper}
p(n;S,M)\le \prod_{a_i\in S}M(n/a_i).
\end{equation}
Further, there must exist at least one integer $r\le n^2$ s.t.
\begin{equation}
\label{eq:lower}
p(r;S,M)\ge
\frac{1}{n^2+1}\prod_{a_i\in S}M(n/a_i).
\end{equation}
If also $p(n;S,M)$ is a nondecreasing function of $n$ then we have
the stronger statement that \[p(n;S,M)\ge \frac{1}{n+1}\prod_{a_i\in
S}M(\sqrt{n}/a_i).\]
\end{lemma}

\noindent Proof. Fix $n>0$ and consider the form
\[\phi=m_1a_1+m_2a_2+m_3a_3+\dots+m_na_n.\]
Now allow each of the $m_i$ to take any value that it wishes to
take, subject to $m_i\in M$ and $m_i\le n/a_i$. For each set of
choices, the form $\phi$ is a partition of some integer $\le n^2$,
and \textit{all} partitions of $n$ occur.

The total number of values that the form takes, counting
multiplicities, is
\[\prod_iM(n/a_i).\]
(Note that all terms with sufficiently large index $i$ are $=1$.)
Since every partition of $n$ occurs, we find that
\begin{equation}
\label{eq:up} p(n;S,M)\le \prod_iM(n/a_i). \end{equation}
Furthermore, since the average number of occurrences of the integers
$\le n^2$ is
\[\frac{1}{n^2+1}\prod_iM(n/a_i),\]
the second conclusion of the lemma is proved. $\Box$

Let's test this with one or two examples. First take $M$ to be all
nonnegative integers and $S$ to be all positive integers. then
$M(x)=1+\lfloor{x}\rfloor$ and we find that $p(n;S,M)\le
\prod_i(1+\lfloor{n/i}\rfloor)$. This is around $n^n/n!$, which is
roughly $e^n$, whereas
 the correct growth is around
$e^{C\sqrt{n}}$. The lower bound is about $e^n/n^2$, so there exists
an integer $r\le n^2$ s.t. $p(r;S,M)\ge
\prod_i\lceil{n/i}\rceil/n^2$, which is about $e^n/n^2$. But indeed,
if there is such an integer $r$, then since $p(n)$ is monotone, we
can take $r=n^2$. Lemma \ref{lem:ineq} then says that $p(n^2;S,M)\ge
e^n/n^2$, or
\begin{equation} \label{eq:lb} p(n;S,M)\ge
e^{\sqrt{n}}/n,
\end{equation} which is reasonably sharp. For another example, in
the case of binary partitions, the upper bound (\ref{eq:up}) yields
the estimate
\[\log{p_S(2n)}\le \log{(2n+1)}{\log_2{(2n)}}\sim
\frac{(\log{n})^2}{\log{2}},\] which can be compared with
(\ref{eq:deb}).

\section{The growth of $p_S(n)$}
\begin{theorem}
Let $S$ be an infinite set of positive integers, and let $p_S(n)$ be the number of
 partitions of $n$ whose parts belong
 to $S$. Then $p_S(n)$ is of superpolynomial growth, that is, for every fixed $k$ the
 assertion $p_S(n)=O(n^k)$ is false.
 This result is best possible in the sense that if $\epsilon(n)$ is any function of
  $n$ that $\to\infty$, then
 we can find an infinite set $S$ such that
 $p_S(n)=O(n^{\epsilon(n)})$.
\end{theorem}
\noindent{Proof.} Let $S=\{1\le a_1<a_2<\dots\}$. Then
$g=$gcd$(S)\le a_1<\infty$, and the theorem is true for $S$ iff it
is true for $S/g$. Hence we can, and do, assume w.l.o.g. that
gcd$(S)=1$.

Let $T\subseteq S$ be such a finite coprime subset, and put $k=|T|$. By Schur's
theorem we have
$p_S(n)\ge p_T(n)\sim Cn^{k-1}$. But we can make $k$ arbitrarily large by adjoining
 elements of $S$
 to $T$ since that adjunction preserves coprimality. Therefore $P_S(n)$ must grow
 superpolynomially.

For the second part of the theorem we use (\ref{eq:upper}) with unconstrained
multiplicities, i.e.,
with $M(x)=1+\lfloor{x}\rfloor$ for $x>0$. If we write $A(n)=|\{i:a_i\le n\}|$ then
 (\ref{eq:upper})
  reads as
\begin{eqnarray*}
p_S(n)&\le& \prod_{i\ge 1}\left(1+\left\lfloor{\frac{n}{a_i}}\right\rfloor\right)\le
\prod_{a_i\le n}\left(1+\frac{n}{a_i}\right)\le n^{A(n)}\prod_{a_i\le n}
\left(\frac{1}{n}+\frac{1}{a_i}\right)
\le n^{A(n)}\prod_{a_i\le n}\left(1+\frac{1}{a_i}\right)\\
&\le&n^{A(n)}\prod_{a_i\le n}e^{1/a_i} \le
n^{A(n)}e^{H_n}=O(n^{A(n)+1}),
\end{eqnarray*}
in which $H_n$ is the $n$th harmonic number. Evidently we can make this
$O(n^{\epsilon(n)})$ by taking the set
$S$ to be sufficiently sparse. $\Box$
\section{A partition function that grows slowly}
There are infinite sequences of allowable parts and multiplicities
on which the partition function grows only polynomially fast, in
fact it can even grow subpolynomially.

One such example is the case where the allowable parts are the
sequence $\{2^{2^j}\}_{j=0}^{\infty}$ and the allowable
multiplicities are \[\{0\}\cup \{2^{2^j}\}_{j=0}^{\infty}.\] In this
case we have, in the notation above,
$M(x)=1+\lfloor{\lg{\lg{x}}}\rfloor,$ for $x\ge 4$, where
``$\lg{}$'' is the log to the base 2. Then by (\ref{eq:up}) we have
\begin{eqnarray*}
p(n;S,M)&&\le \prod_{2^{2^i}\le n/4}
\left(1+\lfloor{\lg{\lg{\frac{n}{2^{2^i}}}}}\rfloor\right)\le
\prod_{2^{2^i}\le
n/4}\left(2\lfloor{\lg{\lg{\frac{n}{2^{2^i}}}}}\rfloor\right)\\
 &&\le (\lg{n})(\lg{\lg{n}})^{\lg{\lg{n}}},
\end{eqnarray*}
which is of sub-polynomial growth. This argument fails if the parts
and multiplicities are all of the powers of 2.

The above argument can be generalized to give a fairly simple
criterion, in terms of the sets of parts and multiplicities, for
polynomial growth of the partition function.
\section{Representing all large integers}
The example above shows that if the allowable multiplicities and
parts are thin enough, even though they both are infinite sets, then
the partition function can grow very slowly. But the example has the
property some arbitrarily large integers are not represented at all.
It may be that if we rule out such situations then the growth must
be superpolynomial. We formulate this as\\
\noindent\textbf{Unsolved problem 1}: Let $S,M$ be infinite sets of
nonnegative integers with $0\notin S$, and let $p(n;S,M)$ be the
number of partitions of $n$ whose parts all lie in $S$ and the
multiplicities of whose parts all lie in $M$. Suppose further that
$p(n;S,M)>0$ for all sufficiently large $n$. Must $p(n;S,M)$ then be
of superpolynomial growth?

\noindent\textbf{Unsolved problem 2}: Find necessary and sufficient
conditions on $S,M$ in order that $p(n;S,M)>0$ for all large enough
$n$. Failing this, find as sharp as possible necessary conditions,
and similarly sufficient conditions for this to happen.

\noindent\textbf{Unsolved problem 3}: Find necessary and sufficient
conditions on $S,M$ in order that $p(n;S,M)$ increase monotonically
for all large enough $n$. Failing this, find as sharp as possible
necessary conditions, and similarly sufficient conditions for this
to happen).

\section{Monotonicity of the partition function}
With reference to unsolved problem 3 above, we consider the case
where the set $S$ of allowable parts is finite and all
multiplicities are allowed, i.e., the problem of Frobenius.

\begin{theorem}
Let $\{p(n)\}$ be generated by
\begin{equation} \label{eq:gen}
G(x)\eqdef\sum_{n\ge 0}p(n)x^n=\frac{1}{\prod_{i=1}^k(1-x^{a_i})},
\end{equation}
 where gcd$(a_1,\dots,a_k)=1$.
The sequence $\{p(n)\}$ is strictly increasing for all sufficiently
large $n$
 if there does not exist a prime $p$ that divides all but one of the
$a_i$'s, i.e., iff every $(k-1)$-subset of the $a_i$'s is coprime.
\end{theorem}

\noindent\textbf{Proof}: Evidently strict monotonicity holds from
some point on iff
\[(1-x)G(x)=\frac{1-x}{\prod_{i=1}^k(1-x^{a_i})}\] has positive power series coefficients, from some
point on. The partial fraction expansion of $(1-x)G(x)$ is of the
form
\begin{eqnarray}(1-x)G(x)&=&\frac{A_0}{(1-x)^{k-1}}+\frac{A_1}{(1-x)^{k-2}}+\dots+\frac{B_0}{(1-\omega
x)^{k_1}}+\frac{B_1}{(1-\omega x)^{k_1-1}}+\dots\\
&&\quad +\frac{C_0}{(1-\zeta x)^{k_2}}+\frac{C_1}{(1-\zeta
x)^{k_2-1}}+\dots.\nonumber
\end{eqnarray}
In the above, $\omega,\zeta,$ etc. run through the primitive $p$th
roots of unity for each prime $p$ that divides one or more of the
$a_i$'s, and $k_1,k_2,\dots$ are the number of $a_i$'s that each of
these primes divides. If no prime divides all but one of the $a_i$'s
then all of the $k_i$'s are $\le k-2$. If in that case we take the
coefficient of $x^n$ on both sides of we have that
\[p(n)-p(n-1)=A_0{n+k-2\choose n}+O(n^{k-3}),\]
which, since $A_0>0$, is positive for all large enough $n$, as
claimed. $\Box$

\section{A refinement of the lower bound}
Let's find a sharper lower bound for $p(n;S)$, when $S$ is an
infinite coprime set of admissible parts, and all multiplicities are
availabe.

Let $A=\{1\le a_1<a_2<\dots<a_k\}$ be a finite coprime subset of
$S$. If we put $r'(n;A)=\sum_{j\le n}p(j;A)$, then an inequality due
to Padberg \cite{pad} states that
\begin{equation}
\label{eq:padb} r'(n;A)\ge \frac{(n+1)^k}{k!a_1\dots a_k}.
\end{equation}

Now, for infinitely many $n$ we have $p(n;A)=\max_{j\le n}p(j;A)$.
Hence for such $n$, $r'(n;A)\le (n+1)p(n;A)$, and therefore
\begin{equation}
\label{eq:ineq} p(n;A)\ge \frac{(n+1)^{k-1}}{k!a_1\dots a_k}.
\end{equation}

Next, extend the set $A$ by adjoining to it the next $h$ basis
elements, to get a new coprime set
\[A_h=\{a_1,a_2,\dots,a_k,a_{k+1},a_{k+2},\dots,a_{k+h}\}.\]
If we apply (\ref{eq:ineq}) to $A_h$ we find that
\[p(n;S)\ge p(n;A_h)\ge \frac{(n+1)^{k+h-1}}{(k+h)!a_1a_2\dots a_{k+h}}.\]

Since $h$ is arbitrary we can optimize this inequality by defining
$j=j(n)$ to be the least integer such that $ja_j\ge n$.
\begin{theorem}
\label{th:lowb}
 Let $S$ be an infinite coprime set, and let $M$ consist of all nonnegative integers.
 Then for large enough $n$ we will have
\begin{equation}\label{eq:shrp}
p(n;S,M)\ge \frac{(n+1)^{j(n)-1}}{(j(n))!a_1a_2\dots a_{j(n)}}.
\end{equation}
\end{theorem}

For example if $S$ consists of all positive integers we find for the
classical partition function that $p(n)\ge e^{2\sqrt{n}}/(2\pi n^2)$
for all large enough $n$, which can be compared to the bound
(\ref{eq:lb}), obtained earlier.


\begin{thebibliography}{aaa}
\bibitem{ber} Bateman, Paul T., and Erd\H os, P\'al, Monotonicity of partition functions, Mathematika, London \textbf{3}, (1956)  1--14.
\bibitem{deb} de Bruijn, N. G., On Mahler's partition problem, Nederl. Akad.Wetensch., Proc. \textbf{51}, (1948) 659--669 = Indagationes Math. \textbf{10}, 210--220 (1948).
\bibitem{hr} Hardy, G. H. and Ramanujan, S, Asymptotic formul\ae\ in combinatory analysis, Proc. London Math. Soc. (2) \textbf{17} (1918), 75--115.
\bibitem{pad} Padberg, Manfred W., A Remark on ``An Inequality for the Number of Lattice Points in a Simplex,'' SIAM J. Appl. Math. \textbf{20} (No. 4) 1971, 638--641.
\bibitem{rad} Rademacher, Hans, On the expansion of the partition function in a series,  Ann. of Math. (2)  \textbf{44},  (1943). 416--422.
    \bibitem{sch} I. J. Schur, \textit{Zur additiven zahlentheorie}, Sitzungsberichte Preussische Akad. Wiss., Phys. Math. Kl. (1926), 488--495.
\end{thebibliography}
\end{document}